\documentclass{amsart}

\usepackage{fancyhdr}
\usepackage{lastpage}
\usepackage{stmaryrd,yhmath}

\newcommand{\bdis}{\begin{displaymath}}
\newcommand{\edis}{\end{displaymath}}
\newcommand{\be}{\begin{equation}}
\newcommand{\ee}{\end{equation}}
\newcommand{\mbb}{\mathbb}
\newcommand{\mcal}{\mathcal}

\newcommand{\vp}{\varphi}

\newcommand{\vth}{\vartheta}

\theoremstyle{definition}

\newtheorem{cor}[]{Corollary}

\theoremstyle{remark}
\newtheorem{remark}[]{Remark}

\newtheorem*{mydef1}{{\bf Theorem}}

\newtheorem*{mydef2}{{\bf Definition}}

\numberwithin{equation}{section}

\begin{document}

\title{Jacob's ladders and the asymptotic formula for short and microscopic parts of the Hardy-Littlewood integral of the function
$|\zeta(1/2+it)|^4$}

\author{Jan Moser}

\address{Department of Mathematical Analysis and Numerical Mathematics, Comenius University, Mlynska Dolina M105, 842 48 Bratislava, SLOVAKIA}

\email{jan.mozer@fmph.uniba.sk}


\begin{abstract}
The elementary geometric properties of Jacob's ladders of the second order lead to a class of new asymptotic formulae for short and microscopic
parts of the Hardy-Littlewood integral of $|\zeta(1/2+it)|^4$. These formulae cannot be obtained by methods of Balasubramanian, Heath-Brown and Ivic.
\end{abstract}

\maketitle

\section{Formulation of the Theorem}

\subsection{}

Let us remind that Hardy and Littlewood started to study the following integral in 1922
\be \label{1.1}
\int_1^T\left|\zeta\left(\frac{1}{2}+it\right)\right|^4{\rm d}t=\int_1^T Z^4(t){\rm d}t,
\ee
where
\bdis
Z(t)=e^{i\vth(t)}\zeta\left(\frac{1}{2}+it\right),\ \vth(t)=-\frac 12 t\ln \pi+\text{Im}\ln\Gamma\left(\frac{1}{4}+i\frac t2\right),
\edis
and they derived the following estimate (see \cite{2}, pp. 41, 59, \cite{14}, p. 124)
\be \label{1.2}
\int_1^T\left|\zeta\left(\frac{1}{2}+it\right)\right|^4{\rm d}t=\mcal{O}(T\ln^4 T) .
\ee
Let us remind furthermore that Ingham, in 1926, has derived the first asymptotic formula
\be \label{1.3}
\int_1^T\left|\zeta\left(\frac{1}{2}+it\right)\right|^4{\rm d}t=\frac{1}{2\pi^2}T\ln^4T+\mcal{O}(T\ln^3T)
\ee
(see \cite{3}, p. 277, \cite{14}, p.125). In 1928 Titchmarsh has discovered a new treatment to the integral (1.1)
\be \label{1.4}
\int_0^T Z^4(t)e^{-\delta t}{\rm d}t\sim\frac{1}{2\pi^2}\frac 1\delta\ln^4\frac 1\delta\Rightarrow \
\int_1^TZ^4(t){\rm d}t\sim\frac{1}{2\pi^2}T\ln^4T
\ee
(see \cite{14}, pp. 136, 143). Let us remind, finally, the Titchmarsh-Atkinson formula (see \cite{14}, p. 145)
\begin{eqnarray} \label{1.5}
& &
\int_0^T Z^4(t)e^{-\delta t}{\rm d}t=\frac 1\delta
\left( A\ln^4\frac 1\delta+B\ln^3\frac 1\delta+C\ln^2\frac 1\delta+\right. \\
& &
\left. +D\ln\frac 1\delta+E\right)+\mcal{O}\left\{\left(\frac 1\delta\right)^{13/14+\epsilon}\right\}, \
A=\frac{1}{2\pi^2} , \nonumber
\end{eqnarray}
which improved the Titchmarsh formula (1.4), and the Ingham - Heath-Brown formula (see \cite{4}, p. 129)
\be \label{1.6}
\int_0^TZ^4(t){\rm d}t=T\sum_{K=0}^4 C_K\ln^{4-K}T+\mcal{O}(T^{7/8+\epsilon}),\ C_0=\frac{1}{2\pi^2}
\ee
which improved the Ingham formula (1.3).

\subsection{}

It is clear that the asymptotic formulae for short and microscopic parts
\be \label{1.7}
\int_T^{T+U}\left|\zeta\left(\frac{1}{2}+it\right)\right|^4{\rm d}t
\ee
of the Hardy-Littlewood integral (1.1) cannot be obtained by methods which lead to the results (1.2)-(1.6). It is proved in this paper that the
Jacob's ladders of the second order $\vp_2(T)$ (see \cite{12}) lead to new asymptotic formulae in this direction. \\

Let us remind our formula (see \cite{12}, (5.11))
\begin{eqnarray} \label{1.8}
& &
Z^4(t)=\frac{1}{2\pi^2}\left\{ 1+\mcal{O}\left(\frac{(\ln\ln T)^2}{\ln T}\right)\right\}\ln^4 T\frac{{\rm d}\vp_2(t)}{{\rm d}t}, \\
& &
t\in [T,T+U_0],\ U_0=T^{13/4+2\epsilon} . \nonumber
\end{eqnarray}
Then from (1.8) the multiplicative asymptotic formula for short and microscopic parts (1.7) of the Hardy-Littlewood integral (1.1) follows
(compare \cite{7}, (1.2)).

\begin{mydef1}
\begin{eqnarray} \label{1.9}
& &
\int_T^{T+U}Z^4(t){\rm d}t=\frac{1}{2\pi^2}\left\{ 1+\mcal{O}\left(\frac{(\ln\ln T)^2}{\ln T}\right)\right\}U\ln^4T
\tan[\alpha_2(T,U)], \\
& &
U\in (0,U_0],\ U_0=T^{13/14+2\epsilon}, \nonumber
\end{eqnarray}
where $\alpha_2$ is the angle of the chord of the curve $y=\vp_2(T)$ that binds the points $[T,\vp_2(T)],\ [T+U,\vp_2(T+U)]$.
\end{mydef1}

\begin{remark}
The small improvements of the exponent $13/14$ of the type $13/14\to 8/9\to \dots$ are irrelevant in this question.
\end{remark}

This paper is a continuation of the series of papers \cite{5}-\cite{13}.

\section{Some canonical equivalences}

\subsection{}

Let us remind that
\be \label{2.1}
\tan[\alpha_2(T,U_0)]=1+\mcal{O}\left(\frac{1}{\ln T}\right)
\ee
is true (see \cite{12}, (5.6)). Then, similarly to \cite{7}, 2.1, we call the chord binding the points
$[T,\vp_2(T)],\ [T+U_0,\vp_2(T+U_0)]$ of the Jacob's ladder $y=\vp_2(T)$ the \emph{fundamental chord} (compare \cite{7}). \\

Let us consider the set of all segments $[M,N]\subset [T,T+U_0]$.

\begin{mydef2}
The chord binding the points
\bdis
[N,\vp_2(N)], [M,\vp_2(M)],\ [M,N]\subset [T,T+U_0],
\edis
such that the property
\be \label{2.2}
\tan[\alpha_2(N,M-N)]=1+o(1),\ T\to\infty
\ee
is fulfilled, is called the \emph{almost parallel chord} to the fundamental chord. This property will be denoted by the symbol $\fatslash_2$,
(comp. \cite{7}).
\end{mydef2}

Now, we obtain the following corollary from (1.9) and (2.2).

\begin{cor}
Let $[M,N]\subset [T,T+U_0]$. Then
\be \label{2.3}
\frac{1}{M-N}\int_N^M Z^4(t){\rm d}t\sim\frac{1}{2\pi^2}\ln^4T\ \Leftrightarrow\ \fatslash_2 .
\ee
\end{cor}

\begin{remark}
We see that the analytic property
\bdis
\frac{1}{M-N}\int_N^M Z^4(t){\rm d}t\sim\frac{1}{2\pi^2}\ln^4T
\edis
is equivalent to the geometric property $\fatslash_2$ of Jacob's ladder $y=\vp_2(T)$ of the second order.
\end{remark}

\subsection{}

Next, similarly to the case of the paper \cite{7}, the following corollary is obtained from our Theorem.

\begin{cor}
There is a continuum of intervals $[M,N]\subset [T,T+U_0]$ such that the asymptotic formula
\be \label{2.4}
\int_N^MZ^4(t){\rm d}t\sim\frac{1}{2\pi^2}(M-N)\ln^4T
\ee
holds true.
\end{cor}

\begin{remark}
Especially, there is a continuum of intervals $[N,M]:\ 0<M-N<1$, such that the asymptotic formula (2.4) is true (this follows from the
elementary mean-value theorem of differentiation).
\end{remark}

\section{On microscopic parts of the Hardy-Littlewood integral (1.1) in neighborhoods of zeroes of the function $\zeta(1/2+iT)$}

Let $\gamma,\gamma^\prime$ be a pair of neighboring zeroes of the function $\zeta(1/2+iT)$. The function $\vp_2(T)$ is necessarily convex
on some right neighborhood of the point $T=\gamma$, and this function is necessarily concave on some left neighborhood of the point $T=\gamma'$.
Therefore, there exists a minimal value $\rho\in (\gamma,\gamma^\prime)$ such that $[\rho,\vp_2(\rho)]$ is the point of inflection of the
curve $y=\vp_2(T)$. At this point, by the properties of the Jacob's ladders, we have $\vp_2^\prime(\rho)>0$. Let furthermore
$\beta=\beta(\gamma,\rho)$ be the angle of the chord binding the points
\be \label{3.1}
[\gamma,\vp_2(\gamma)],\ [\rho,\vp_2(\rho)] .
\ee
Then we obtain by Theorem (compare \cite{7})
\begin{cor}
For every sufficiently big zero $T=\gamma$ of the function $\zeta(1/2+iT)$ the following formulae describing microscopic parts (1.7) of the
Hardy-Littlewood integral (1.1) hold true
\begin{itemize}
\item[(A)] a continuum of asymptotic formulae
\begin{eqnarray} \label{3.2}
& &
\int_\gamma^{\gamma+U}Z^4(t){\rm d}t\sim \frac{\tan\alpha}{2\pi^2}U\ln^4\gamma,\ \gamma\to\infty, \\
& &
\alpha\in (0,\beta(\gamma,\rho)),\ U=U(\gamma,\alpha)\in (0,\rho-\gamma) , \nonumber
\end{eqnarray}
where $\alpha=\alpha(\gamma,U)$ is the angle of the rotating chord binding the points $[\gamma,\vp_2(\gamma)],\ [\gamma+U,\vp_2(\gamma+U)]$,
\item[(B)] a continuum of asymptotic formulae for a chord parallel to the chord given by the points (3.1)
\be \label{3.3}
\int_N^MZ^4(t){\rm d}t\sim\frac{\tan[\beta(\gamma,\rho)]}{2\pi^2}(M-N)\ln^4\gamma,\ \gamma<N<M<\rho .
\ee
\end{itemize}
\end{cor}

\begin{remark}
Let us remind that if the Riemann conjecture  is true then the Littlewood estimate
\bdis
\gamma^\prime-\gamma<\frac{A}{\ln\ln\gamma}\to 0,\ \gamma\to\infty
\edis
takes place (a simple consequence of the estimate $S(T)=\mcal{O}(\ln T/\ln\ln T)$, see \cite{14}, p. 296).
\end{remark}

\section{Second class of formulae for parts (1.7) of the Hardy-Littlewood integral (1.1) beginning in zeroes of the function $\zeta(1/2+iT)$}

Let $T=\gamma,\bar{\gamma}$ be a pair of zeroes of the function $\zeta(1/2+iT)$, where $\bar{\gamma}$ obeys the following conditions
(compare \cite{7})
\bdis
\bar{\gamma}=\gamma+\gamma^{13/14+2\epsilon}+\Delta(\gamma),\ 0\leq \Delta(\gamma)=\mcal{O}(\gamma^{1/4+\epsilon})
\edis
(see the Hardy-Littlewood estimate for the distance between the neighboring zeroes \cite{1}, pp. 125, 177-184). Consequently
\be \label{4.1}
U(\gamma)=\gamma^{13/14+2\epsilon}+\Delta(\gamma)\sim ^{13/14+2\epsilon},\ \gamma\to\infty .
\ee
For the chord that binds the points
\be \label{4.2}
[\gamma,\vp_2(\gamma)],\ [\bar{\gamma},\vp_2(\bar{\gamma})]
\ee
we obtain (similarly to \cite{12}, (5.6))
\be \label{4.3}
\tan[\alpha_2(T,U)]=1+\mcal{O}\left(\frac{1}{\ln T}\right) .
\ee
The continuous curve $y=\vp_2(T)$ lies below the chord given by the points (4.2) on some right neighborhood of the point $T=\gamma$, and this
curve lies above that chord on some left neighborhood of the point $T=\bar{\gamma}$. Therefore there exists a common point of the curve and the
chord. Let $\bar{\rho}\in (\gamma,\bar{\gamma})$ be such a common point that is the closest one to the point $[\gamma,\vp_2(\gamma)]$. Then we obtain
from our Theorem (compare \cite{7}) the next corollary.
\begin{cor}
For every sufficiently big zero $T=\gamma$ of the function $\zeta(1/2+iT)$ we have the following formulae for the parts (1.7) of the
Hardy-Littlewood integral (1.1)
\begin{itemize}
\item[(A)] a continuum of asymptotic formulae for the rotating chord
\begin{eqnarray} \label{4.4}
& &
\int_\gamma^{\gamma+U}Z^4(t){\rm d}t\sim\frac{\tan\alpha}{2\pi^2}U\ln^4\gamma,\ \tan\alpha\in [\eta,1-\eta], \\
& &
U=U(\gamma,\alpha)\in (0,\bar{\rho}-\gamma) , \nonumber
\end{eqnarray}
where $\alpha$ is the angle of the rotating chord binding the points $[\gamma,\vp_2(\gamma)]$ and $[\gamma+U,\vp_2(\gamma+U)]$, and
$0<\eta$ is an arbitrary small number,
\item[(B)] a continuum of asymptotic formulae for the chords parallel to the chord binding the points (4.2), (see (4.3))
\be \label{4.5}
\int_N^M Z^4(t){\rm d}t\sim\frac{1}{2\pi^2}(M-N)\ln^4\gamma,\ \gamma\leq N<M\leq \bar{\rho} .
\ee
\end{itemize}
\end{cor}

\begin{remark}
For example, in the case $\alpha=\pi/6$ we have from (4.4)
\bdis
\int_\gamma^{\gamma+U} Z^4(t){\rm d}t\sim \frac{1}{2\sqrt{3}\pi^2}U\ln^4\gamma,\ U=U\left(\gamma,\frac{\pi}{6}\right) .
\edis
\end{remark}

\begin{remark}
It is obvious that (see (4.4))
\bdis
U(\gamma,\alpha)<T^{7/8+2\epsilon} .
\edis
Moreover, the following is also true
\bdis
U(\gamma,\alpha)<T^{\omega+2\epsilon},\ \omega<\frac{7}{8} ,
\edis
where $\omega$ is an arbitrary improvement of the exponent $7/8$ which will be proved.
\end{remark}

\begin{remark}
The asymptotic formulae (1.9), (2.3), (2.4), (3.2), (3.3), (4.4), (4.5) cannot be derived within complicated methods of
Balasubramanian, Heath-Brown and Ivic (compare \cite{4}).
\end{remark}

\thanks{I would like to thank Michal Demetrian for helping me with the electronic version of this work.}

\end{document}